\documentclass[11pt]{amsart}
\usepackage{amsmath}
\usepackage{amsfonts}
\usepackage{amssymb}
\usepackage{amsthm}
\usepackage{color}

\newcommand{\CC}{{\mathbb C}}
\newcommand{\FF}{{\mathbb F}}
\newcommand{\ZZ}{{\mathbb Z}}
\newcommand{\Fq}{{\FF_q}}
\newcommand{\unit}{{\bf 1}}
\newcommand{\Cent}{{\mathrm C}}
\newcommand{\compl}{{\rm c}}
\newcommand{\GL}{{\rm GL}}
\newcommand{\PGL}{{\rm PGL}}
\newcommand{\sgnCD}{{{\rm sgn}_{\rm CD}}}
\newcommand{\bGL}{{\bf GL}}
\newcommand{\Id}{{\rm I}}
\newcommand{\Ind}{{\rm Ind}}
\newcommand{\Inv}{{\rm Inv}}
\newcommand{\inv}{{\rm inv}}
\newcommand{\Pair}{{\rm Pair}}
\newcommand{\sign}{{\rm sign}}
\newcommand{\sgn}{{\rm sgn}}
\theoremstyle{plain}
\newtheorem{thm}{Theorem}
\newtheorem{prop}{Proposition}[section]
\newtheorem{prop3}{Proposition}[subsection]

\newtheorem{lem}[prop]{Lemma}

\theoremstyle{definition}
\newtheorem{ex}[prop]{Example}

\newtheorem{dfn}[prop]{Definition}
\newtheorem{rmk}[prop]{Remark}

\newcommand\rC{{\mathrm C}}
\newcommand\rD{{\mathrm D}}
\newcommand\rE{{\mathrm E}}
\newcommand\rG{{\mathrm G}}
\newcommand\rF{{\mathrm F}}
\newcommand\rH{{\mathrm H}}
\newcommand\rI{{\mathrm I}}
\newcommand\Green{{\mathrm{Green}}}
\newcommand\Isotr{{\mathrm{Isotr}}}
\newcommand\SL{{\mathrm{SL}}}

\begin{document}

\title{ Groupoids, Geometric Induction and Gelfand Models} 
\author{Anne-Marie Aubert } 
\address{Sorbonne Universit\'e and Universit\'e de Paris, CNRS,
Institut de Math\'ematiques de Jussieu -- Paris Rive Gauche, IMJ-PRG, F-75006, Paris, France}
\email{anne-marie.aubert@imj-prg.fr}
\author{Antonio Behn} 
\address{  Facultad de Matem\'aticas, Pontificia Universidad Cat\'olica de Chile, Santiago, Chile}
\email{antonio.behn@mat.puc.cl}
\thanks{The third author was partially supported by Fondecyt Grant 1140510}
\author{Jorge Soto-Andrade} 
\address{Departamento de Matem\'aticas, Facultad de Ciencias, Universidad de Chile, Santiago, Chile}
\email{sotoandrade@uchile.cl}
\date{\today}



\begin{abstract}
In this paper we introduce an intrinsic version of the classical induction 
of representations for a subgroup $H$ of a (finite) group $G$, called here 
{\em geometric induction}, which associates to any, not necessarily transitive,  
$G$-set $X$ and any representation of the action groupoid $A(G,X)$ associated 
to $G$ and $X$, a representation of the group $G$.  
We show that geometric induction, applied to one dimensional characters of the 
action groupoid of a suitable $G$-set $X$ affords a Gelfand Model for $G$ in 
the case where  $G$ is 
either the symmetric group or the projective general linear group of rank $2$. 

\end{abstract}

\maketitle


\section{Introduction}

A groupoid is a natural extension of a group, considered as a
category: it is a category where each morphism is invertible. A groupoid  with only one object is a group.

Mackey's theory of induced representations from a subgroup $H$ of a group $G$ 
to $G$ itself, may be described in intrinsic, more geometric terms. The idea 
is to associate a representation of $G$ not to a representation of a subgroup 
$H$ of $G$ but to a $G$-set $X$ and a representation $\sigma_{G,X}$ of the so called 
 {\em action groupoid} $A(G,X)$, also denoted
$X\rtimes G$, associated to $X$.             

Recall that  the one dimensional characters of  $A(G,X)$  play a key role in the geometric construction   of  linear representations of   $G$. 
Indeed, we can twist the natural representation of $G$ in the function space  
$L^2(X)$ by those characters, to obtain multiplicity-free representations of  
$G$ in many cases, even if the original action is not transitive.

We show below that  in this way we can obtain  a geometric construction of 
Gel'fand Models for symmetric groups and the projective general linear group of rank $2$. 

On the other hand, groupoids also play a significant role in group representation
theory when we look at representations of a whole family of groups instead of 
just one. For instance, when  we construct linear representations of  the 
``field'' of symmetric groups $S_n$, each sitting on a set of size $n$, for 
all $n$. More precisely, we consider then the groupoid $\bf S$ whose objects are 
finite sets, and whose arrows are  bijections among finite sets. We have a 
$q$-analogue of this, i.e. the linear groupoid  $\bGL(q)$ 
whose objects are finite dimensional vector spaces over $\mathbb F_q$ and whose 
arrows are linear isomorphisms.
 

\section{Induced representations: a geometric approach}
Let $G$ be a finite group. If $G$ acts on a finite set $X$, we can construct the 
associated action groupoid  $ A(G,X) $ in  the following way:

\begin{itemize}
\item[--]
the objects of  $ A(G,X) $ are the points of  $X$.
\item[--]
the arrows of  $ A(G,X) $  are the triples  $(x, g, y)$ such that   
$y = g\cdot x$, with the obvious composition. 
\end{itemize}

Given an action  of a group $G$ on a set $X$, the action groupoid $A(G,X)=
X\rtimes G$ is closely related to the quotient set $X/G$ (the set of $G$-orbits). 
But, instead of
taking elements of $X$ in the same $G$-orbit as being equal in $X/G$, in the
action groupoid they are just isomorphic.

 
\smallskip
 
A {\em representation} of the groupoid $A(G,X)$ is an action of $G$ on a finite dimensional vector bundle $p\colon W \to X$. 
A representation gives a functor $\sigma$ from $A(G,X)$ to the category of finite dimensional complex vector spaces such that $\sigma(x)$ is the fiber $W_x$ of $p$ at $x$ and $\sigma(x,g,y)\colon W_x\to W_y$ is an isomorphism for each $g\in G$ such that $g\cdot x=y$.

Now, to any representation  $(\sigma, W)$ of   $ A(G,X) $  we can associate 
the induced representation $(\rho,V) $ of  $\sigma$ to $G$, which we call {\em the natural representation $\rho$ of  $G$ associated to  the  $G$-set  $X$ twisted by  $\sigma$},  denoted   $\Ind_{A(X,G)}^G(\sigma) $,   defined as follows:
\begin{itemize}
\item[(a)]
the space of  $\rho$ is  the complex vector space   $ V =  L^{2}_{W}(X)$ of all sections of $p$; 
\item[(b)]
$ (\rho \sb{g} f )(x) =  \sigma(g^{-1}\cdot x, g, x) [ f(g^{-1}\cdot x)]$          
for all  $ g \in G$,  $x \in X$,  $f \in  V$. 
\end{itemize}

\begin{prop}
With the above notations, if   $X = \bigcup_{i \in I} X_{i} $ is the  $G$-orbit decomposition of   $X$ and $H_{i}$ is the isotropy group of a chosen point  
$x_{i} \in X_{i}$  for each  $i \in I $, then we have
\[(\rho,V) \cong     \bigoplus_{i \in I} \Ind_{H_i}^G(\sigma|_{H_i}),\]
where   $\sigma |_{H_{i}}$ stands for the restriction of the representation $\sigma$ to the 
isotropy group $H_{i}$, defined by    
$\sigma |_{H_{i}} (h) := \sigma  (x_{i}, h, x_{i})$  for  $  h \in H_{i}$.      
\end{prop}
   \begin{flushright}    
    $ \square   $  
\end{flushright} 

Notice that it follows from this proposition that the induced representation $(\rho,V)$ depends only on the restriction of the inducing character $\sigma$ to the isotropy groups $H_i.$

\section{Gelfand models}

A Gel'fand model for a finite group $G$ is a complex linear representation $M$ of $G$ that contains each of its irreducible representations with multiplicity
one. The representation $M$ is called an {\em involution model} if there exists a set of representatives 
$\{\omega_i\}$ of the distinct conjugacy classes of involutions in $G$ and a set 
$\{\lambda_i \colon \rC_G(\omega_i) \to\CC \}$ of linear characters of
the centralizers $\rC_G(\omega_i)$ in $G$ of the $\omega_i$, such that   
$$ M \simeq \bigoplus_i \Ind_{\rC_G(\omega_i)}^G(\lambda_i)$$ 

Klyachko \cite{K}, Inglis, Richardson and Saxl \cite{IRS}, Kodiyalam and Verma
\cite{KV}, Adin, Postnikov and Roichman \cite{APR} constructed {\em 
involution models} for the symmetric group. 
Baddeley classified which irreducible Weyl groups admit involution models:
in particular only those of type $\rD_{2n}$ ($n>1$);  $\rF_4$, $\rE_6$, $\rE_7$
and $\rE_8$ do not have involution models. Vinroot extended this
classification to finite Coxeter groups, in particular showing that the
Coxeter groups of type $\rI_2(n)=\rG_2(n)$ (the dihedral groups) and $\rH_3$ do
have involution models while the group $\rH_4$ does not.

Bump and Ginzburg introduced in \cite{BG} the more flexible notion of {\em generalized involution
model}.  
Marberg proved \cite{M2} that a finite complex reflection group has a
generalized involution model if and only if each of its irreducible
factors is one of the following: 
$\rG(r,p,n)$ with $\gcd(p,n)=1$, $\rG(r,p,2)$ with $r/p$ odd, or $\rG_{23}$, the
Coxeter group of type $\rH_3$.

On the other hand, Garge and Oesterl\'e constructed in \cite{GO}
for any finite group $G$ with a faithful representation $V$, 
a representation which they called the {\em polynomial model for $G$ 
associated to $V$}. Araujo and others have proved that the polynomial
models for certain
irreducible Weyl groups associated to their canonical representations are
Gel'fand models.
In \cite{GO} it is proved that a polynomial model
for a finite Coxeter group $G$ is a Gel'fand model if and only if $G$ has no
direct factor of the type $\rD_{2n}$, $\rE_7$ or $\rE_8$.

Combining the above results of \cite{M2} and \cite{GO} shows that 
\begin{itemize}
\item[-]
the groups $\rF_4$, $\rH_4$ and $\rE_6$ (which in the Shephard-Todd
classification are $\rG_{28}$, $\rG_{30}$ and $\rG_{35}$, respectively) do not have a
generalized involution model, however they have a polynomial model which is a
Gel'fand model;
\item[-]
the groups $\rE_7$ and $\rE_8$ (which in the Shephard-Todd
classification are $\rG_{36}$ and $\rG_{37}$, respectively) do not have a
generalized involution model, have a polynomial model, but this polynomial
model is not a Gel'fand model.
\end{itemize}

\subsection{A geometric Gelfand model for the symmetric group}  
Let  $G=S_n$ be the symmetric group realized as the permutation group of $\{ 1, 2, \ldots, n\}$ and 
take  for $X$ the set of all involutions of $G$.

In order to illustrate the general picture, we start by the example of the
group $S_4$.                                                                                                                      
\begin{ex}
Let $G=S_4$. Let $\Id$ denote the identity in $G$. We have
\[X=\left\{\Id, (12), (13), (14), (23), (24), (34), (12)(34), (13)(24), (14)(23)\right\}.\] 
Then   $|X|=10$ and $G$ has three orbits in $X$,  to wit, the orbit that is 
reduced to the identity $\Id$, the orbit consisting of all $6$ transpositions 
$(ij) $ and  the orbit  consisting of all $3$ double transpositions  $(ij)(kl)$.  
Let  $\sigma=\sigma_{X,S_4}$ denote the one dimensional representation of $A(X,G) $ defined by 
\[\sigma(x, g) :=\begin{cases} 
1&\text{if $x=\Id$, for all $ g  \in  G$;}\cr
1&\text{if $x=(ij)$ with $i<j$ and $g(i)<g(j)$;}\cr   
-1&\text{if $x=(ij)$  with $i<j$ and $g(i)>g(j)$;}\cr
1&\text{if $x=(ij)(kl)$ with $i<j,k<l$ and $g(i)<g(j)$ and $g(k)<g(l)$ or  
$g(i)>g(j)$;}\cr
&\text{and  $g(k)>g(l)$} \cr
-1&\text {otherwise.}
\end{cases}\] 
  Then the isotropy group   $G_{\Id}$  at  $\Id$ is $G$,  the isotropy group 
$G_{(12)}$ at  $(12)$ is  
 $\left\{\Id,(12),(34),(12)(34)\right\} \simeq C_{2}\times C_{2}$, where
$C_2\simeq\ZZ/2\ZZ$ is the cyclic group of order $2$, and 
  the isotropy group $G_{(12)(34)}$ at $(12)(34)$ is  
$$\left\{\Id,(12),(34),(12)(34),(14)(23),(13)(24),(1324),(1423)\right\}\,
\simeq\,(C_{2}\times C_{2}) \rtimes C_{2}.$$
  Notice that these isotropy groups appear as automorphisms groups of the trees  naturally associated to the given partitions.

 The restriction of  the representation   $\sigma_{X,S_4}$   to   $G_{\Id}$ is the unit character  ${\bf 1}$  of  $G$, the restriction of  $\sigma_{X,S_4}$ to   $G_{(12)}  $  is the tensor product  $\sgn \otimes \bf 1 $ of the non trivial character $\sgn$ of  
$C_{2}$ and the trivial character  $\bf 1$ of  $C_{2}$, the restriction of  
$\sigma_{X,S_4}$ to $G_{(12)(34)}$ is the restriction of the signature character $\sign$ of $G$ to 
$G_{(12)(34)}$.
 So we obtain that the representation $(\rho,V)$ of $G$ induced from $\sigma_{X,S_4}$ is 
isomorphic to
$${\bf 1} \oplus  \Ind_{C_{2}\times C_{2}}^{ G} (\sgn \otimes {\bf 1}) \oplus  
\Ind_{(C_{2}\times C_{2})\rtimes C_{2}}^{ G} (\sign).  $$

Recall now that $G$ has $5$ irreducible representations, to wit:
\begin{itemize}
\item[-]
the unit representation $\unit$;
\item[-]
the signature representation  $\sign$;
\item[-]
a $2$-dimensional representation, that we denote   by   $\sigma^{(2)}$;
\item[-]
a $3$-dimensional representation, denoted   $\sigma^{(3)} $, such that   $ {\bf 1} \oplus \sigma^{(3)} \simeq
\tau $, where  $\tau$ stands for the natural (permutation) representation of  
$G$ associated to its canonical action on a $4$-point set;
\item[-]
the $3$-dimensional representation,    $\sign \otimes \sigma^{(3)} $. 
\end{itemize}

It is easy to check that:
$$\Ind_{C_{2}\times C_{2}}^{G} (\sgn \otimes 1)  \simeq  \sigma^{(3)}   \oplus   (\sign \otimes \sigma^{(3)} )\quad\text{and}\quad
\Ind_{(C_{2}\times C_{2})\rtimes C_2}^{G} (\sign )  \simeq   
\sign \oplus  \sigma^{(2)}.$$
It follows that the induced representation $(\rho,V)$ is indeed a Gel'fand 
Model for $G=S_4$.
\end{ex}

\smallskip

We will now consider the case of $G=S_n$ with $n$ arbitrary.

\smallskip

For any positive integer $t$, we will denote by $C_2\wr S_t$ 
the {\em wreath product} of $C_2$ by $S_t$, that is,
the semidirect product $(C_2)^t\rtimes S_t$, where the symmetric group $S_t$
acts on $(C_2)^t$ by permuting the $t$ factors $C_2$. 
In other words
$C_2\wr S_t$ is the {\em hyperoctahedral group} $W(B_t)$ (the Weyl group of type
$B_t$).

For $0\le t\le \lfloor n/2 \rfloor$, let $X_t$ denote the conjugacy class of
involutions of $G$ that are product of $t$ disjoint
transpositions. We have $$X=\bigsqcup_{t=0}^{\lfloor n/2 \rfloor}X_t.$$
Let $x\in X_t$, with $0\le t\le\lfloor n/2 \rfloor$. The involution $x$ fixes 
precisely $f=n-2t$ points. 

We write in the cycle notation:
$$x=(i_1j_1)(i_2j_2)\cdots(i_{t}j_{t}),$$ where $i_{k}<j_{k}$ for
each $k\in\{1,2,\ldots,t\}$. Then the isotropy group $G_x=\Cent_G(x)$ 
at $x$ is generated 
by permutations of following forms:
\begin{itemize}
\item[-] $(i_ki_{k+1})(j_kj_{k+1})$, for $1\le k\le t$;
\item[-] $(i_kj_k)$, for $1\le k\le t$;
\item[-] any permutation of $S_n$ which
fixes the $2t$ points $i_1$, $j_1$, $\ldots$, $i_t$, $j_t$, in other
words,
any permutation in $S_f$.
\end{itemize} 
It follows that $G_x$ is isomorphic to the group
$W(B_t)\times S_f$, where $W(B_t)=C_2\wr S_t$ is embedded in $S_n$ so that 
the subgroup $(C_2)^t$ is generated by the $2$-cycles of $x$.
 
As in \cite{APRII} and \cite{M}, we attach to each permutation  $g\in G$ 
its {\em inversion set} $\Inv(g)$ and
the set $\Pair(g)$ of the $2$-cycles in $g$, that is:
$$\Inv(g):=\left\{(i,j)\,:\,1\le i<j\le n,\;g(i)>g(j)\right\},$$ 
$$\Pair(g):=\left\{(i,j)\,:\,1\le i<j\le n,\;g(i)=j,\;g(j)=i\right\}.$$

Let $x\in X$ and $g\in G$.
We set
$$\Inv_x(g):=\Inv(g)\cap \Pair(x),$$
and we will denote by $\inv_x(g)$ the cardinality of the set $\Inv_x(g)$.

\begin{dfn} Let $\sigma_{G,X}\colon A(G,X)\to \CC$ be the map defined by
$$\sigma_{G,X}(x,g):=(-1)^{\inv_x(g)},\quad\text{for $x\in X$ and $g\in G$.}$$
\end{dfn}

\begin{rmk}
Let $x$ be a involution in $S_n$ which is the product of $t$ disjoint
transpositions. If we write $x$ as
$x=(i_1i_2)(i_3i_4)\cdots(i_{2t-1}i_{2t})$, where $i_{2j-1}<i_{2j}$ for
each $j=1,2,\ldots,t$.  Then $\inv_x(g)$ equals the cardinality of
$$\left\{j\in [1,t]:g(i_{2j-1})>g(i_{2j})\right\},$$
for any $g\in G$. Hence the sign which occurs is the definition of
$\sigma_{G,X}$ coincides with the sign which was introduced by Kodiyalam and 
Verma in \cite{KV}. It is also the sign which occurs in \cite{APR} and \cite{M}.
\end{rmk}

\begin{lem} The map $\sigma_{G,X}$ defines a one dimensional
representation of $A(X,G)$.
\end{lem}
\begin{proof}
Let $x_1,x_2\in X$ and let $g_1,g_2\in G$ such that $x_1=g_2\cdot x_2$.
Let $P_n:=\left\{(i,j)\,:\,1\le i<j\le n\right\}$. For any subset $Q$ of
$P_n$, we will denote by $Q^{\compl}$ its complement in $P_n$. 
Let $i$ and $j$ two elements of $\{1,\ldots,n\}$. We have 
$$(i,j)\in \Inv_{x_2}(g_1g_2)\cap \Inv(g_2)^{\compl}\quad\text{
if and only if }\quad
\text{$i<j$ and $(g_2(i),g_2(j))\in \Inv_{x_1}(g_1)$;}\leqno{-}$$
$$(j,i)\in \Inv_{x_2}(g_2)\cap\Inv(g_1g_2)^\compl\quad\text{
if and only if }\quad
\text{$i>j$ and $(g_2(i),g_2(j))\in \Inv_{x_1}(g_1)$.}\leqno{-}$$
It follows that
$$\inv_{x_2}(g_1g_2)\equiv \inv_{x_2}(g_2)\,+\,\inv_{x_1}(g_1)
\pmod 2.$$
We obtain
$$\sigma_{G,X}(x_2,g_2)\circ\sigma_{G,X}(x_1,g_1)=\sigma_{G,X}(x_2,g_1g_2),$$
that is, $\sigma_{G,X}$ is a representation of $A(X,G)$.
\end{proof}

There is a unique character $\sgnCD\colon W(B_t)\to\{\pm 1\}$ 
whose restriction to the
normal subgroup $(C_2)^t$ of $W(B_t)$ is the product of the sign 
characters of $C_2$ and
that is trivial on the subgroup $S_t$. The kernel of $\sgnCD$ is
isomorphic to the Weyl group $W(D_t)$.
In the parametrization of the irreducible representations of $W(B_t)$ by
the pairs of partitions of $t$ (see \cite{L}) the representation
afforded by the character $\sgnCD$ corresponds to $(\emptyset,(t))$.

Let $x\in X$. Recall that $$G_x=W(B_t)\times S_f.$$
The restriction of $\sigma_{G,X}$ to $G_x$ is the tensor product of
$\sgnCD$ by the unit representation of $S_f$:
$$(\sigma_{G,X})_{|G_x}=\sgnCD\otimes \unit_{S_f}.$$
By transitivity of induction, we get
$$\Ind_{G_x}^G (\sigma_{G,X})_{|G_x}=
\Ind_{S_{2t}\times S_f}^{S_n}(\Ind_{G_x}^{S_{2t}\times S_f}
(\sgnCD\otimes \unit_{S_f}))=\Ind_{S_{2t}\times S_f}^{S_n}\left(
\Ind_{W(B_t)}^{S_{2t}}(\sgnCD)\,\otimes \unit_{S_f}\right).$$
The induced representation $\Ind_{W(B_t)}^{S_{2t}}(\sgnCD)$ is the
mutiplicity free sum of the irreducible characters of $S_{2t}$
corresponding to partitions of $2t$ with even columns only (see, for
instance, \cite[Ch.I \S8 and Ch. VII(2.4)]{McD}).
By using the Littlewood-Richardson rule, it follows that
$\Ind_{G_x}^G \big(\sigma_{G,X})_{|G_x}\big)$ is the multiplicity free sum of all
irreducible Specht modules indexed by partitions with exactly $f$ odd columns. 
It follows that the induced representation $(V, \rho)$ is indeed a Gel'fand 
Model for $G=S_n$.
 
\subsection{A geometric Gel'fand Model for the group  $G = \PGL_2(q)$,  $q$  
odd}
In  \cite{BG} a Gelfand Model for  $G = \PGL(2,q)$ is constructed as a sum of 
three induced representations, taking advantage of the results of \cite{sa:a}. 
More precisely, we have
the natural action of $G$ on the set $X$ consisting of all symmetric matrices 
in  $G$, which may be identified with non degenerate symmetric bilinear forms mod centre or non-degenerate quadratic forms mod centre on the finite plane over  $ k = {\mathbb F}_q  $.   This action has 3 orbits:

\begin{itemize}
\item[-]  the orbit $O_1$ of the split isotropic hyperbolic form $H$ mod centre, 
given by  $ (x,y) \mapsto xy $;
\item[-]  the orbit  $O_2$ of the non-split anisotropic normic form mod centre 
given by the norm $N$ of the unique quadratic extension of  $k$.
\item[-]
the orbit  $O_3$ of the alternating form $\det$ mod centre. 
\end{itemize}

The isotropy groups  $H_1$, $H_2$ and $H_3$ in $G$ for these orbits are just the 
corresponding two projective similarity orthogonal groups and the projective 
similarity symplectic group; more precisely:

\begin{itemize}
\item[-]
$H_1$  is the normalizer of the split torus  $T_1$ in  $G$;
\item[-]
$H_2$ is the normalizer of the non-split torus  $T_2$ in  $G$;
\item[-]
$H_3$ is the whole group  $G$. 
\end{itemize}

Bump and Ginzburg make an {\em ad hoc} choice of characters  $ \psi_1$, 
$\psi_2 $ and $\psi_3$ to induce from  $H_1,   H_2 $ and   $H_3$: 

-   $ \psi_1$ is the trivial character  {\bf 1}  for  $H_1$;

-  $ \psi_2 $ is the order $2$ character of  $H_2$ that is trivial on the 
non-split torus

- $ \psi_3 $   is  the non unique trivial linear character of  $G$, i. e. the ``  sign  character''  given by the determinant mod squares.

\begin{thm}
A model $M$ for  $G$ is given by

$$  M = \Ind_{H_1}^G \psi_1  \oplus    \Ind_{H_2}^G \psi_2 \oplus \Ind_{H_3}^G \psi_3 $$ 

\end{thm}

 \begin{flushright}    
    $ \square   $  
\end{flushright} 
  
\indent 

Our viewpoint is to obtain  the model  $M$ as a geometric induced representation from a `` sign character'' $\varepsilon $ of the action groupoid   $A(X, G)$ naturally associated to the  $G$-space $X$.    

The sign character  $\varepsilon $, which takes values  $ \pm 1$ only,  is defined first on the holonomy (isotropy) groups  $G_x$ associated to any  quadratic form   $x \in X$ as follows:

 \begin{dfn}  Let $\Isotr(x) $ denote the set of all isotropic $1$-dimensional vector subspaces of $x \in X$. Then  $G_x$ acts naturally on  $\Isotr(x)$ and we define   
\[\varepsilon(g) :=\text{sign of the permutation of  $\Isotr(x)$ defined by  $g \in  G_x$.}\]
  \end{dfn}
  
Denote by    $  \varepsilon_i $ the restriction of  $\varepsilon$ to  $H_i$, for $ i = 1, 2, 3.$ Then
  \begin{itemize}
	\item[$\bullet$]
  $\varepsilon_1$ is the order $2$ character   of $H_1$ which is trivial on $T_1$;     
  \item[$\bullet$]
  $\varepsilon_2$ is  the trivial character  of  $H_2$;
  \item[$\bullet$]
  $\varepsilon_3 $  is the  sign character   of   $G$.
  \end{itemize}
	
  Indeed,  the determinant mod squares character at   $g \in G$ may be obtained as well as the signature of the permutation induced by  $g$   on the set of all $1$-dimensional subspaces of the finite plane (the isotropic lines for det mod centre). 
  
  Then for any extension $\tilde \varepsilon $ of the character $\varepsilon$ to the whole action groupoid  $A(X, G),$   we have

  \begin{thm}
  The Model $M$ may be obtained   as
  \[M =  \Ind_{A(X,G)}^G(\tilde \varepsilon)          \simeq \Ind_{H_1}^G \varepsilon_1  \oplus    \Ind_{H_2}^G \varepsilon_2 \oplus \Ind_{H_3}^G \varepsilon_3.\]  
  \end{thm}
  
  $\hfill \square$
  
  Recall that  the induced representation $\Ind_{A(X,G)}^G(\tilde \varepsilon) $
  depends only on $\varepsilon.$ 
  
  Notice also that we may extend  $\varepsilon$ to a full-fledged character of the action grupoid $A(X, G)$, by choosing first a point $x_i$ in each orbit $O_i$, then  for each $x \in O_i, \; x \neq x_i$ an element $g_x \in G$ such that  $g_x \cdot x_i = x$ and assigning  arbitrarily  a non-zero value $\alpha_x$ to each arrow $(x_i, g_x, x)$ in $A(X, G)$.
\subsection{Gelfand-Graev representations as geometrically induced representations: The case of  $\GL_2(q)$}
To construct the classical Gel'fand-Graev representation of $G=\GL_2(q)=\GL_2(\Fq)$ as a 
geometrically induced representation, we consider first the  $G$-set $X$ 
consisting of all pairs $(u,\omega)$,  where  $u$ is a non zero vector in $k^{2}$
and $\omega$ is a non zero $2$-vector in  $\bigwedge^{2}k^{2}$. The elements of  
$X$ may be called Grassmann chains. The group $G$ acts in a natural way on $X$. 

So we have a  $G$-set  $X$ and the associated action groupoid  $M = A(X, G)$.

We proceed now to define a linear character $\widetilde  \psi $ of the groupoid  
$M$ associated to any character $\psi$ of the additive group of the finite field
$k=\mathbb F_{q}$.
Let  $\psi$  be a non trivial (complex) character of the additive group of the finite field  $k = \mathbb F_{q}$.
Choose arbitrarily, for each  $(u,\omega) \in  X$, a vector $v \in k^{2}$, such 
that  $ \omega = u \wedge v$.

Notice that $v$ is well defined modulo  the line  $<u>$  generated by $u.$ So  we are choosing in fact an element in   the fiber $  k^2 /<u> $  above   $u$. 
We have then a natural   $G-$fiber bundle  $ \mathcal E $ with the set of all lines through the origin in   the finite plane   $k^2$  as basis.  

Then, for  $(x, g, y) \in M $, if we write  $ x = (u, \omega) = (u, u \wedge v)$,
$y = (u', \omega' )= (u', u' \wedge  v' )$, there is a unique $b \in k $ such 
that  $g(v) = v' +  bu' $.  
Define   
$$ \widetilde  \psi (x,g,y) = \psi(b) $$
with $b$ as above.
\begin{prop3}
The mapping $\widetilde  \psi $ just defined is a linear character of the 
groupoid $M$.

\end{prop3}
\begin{proof}
We have to check that 
$$  \widetilde  \psi  ((x,g,y)(y,h,z)) =  \widetilde  \psi  (x,g,y) \cdot  \widetilde  \psi  (y,h,z), $$   i. e. 
$$      \widetilde  \psi  ((x, h \circ g,z)) =  
\widetilde  \psi  (x,g,y) \cdot  \widetilde  \psi  (y,h,z). $$
Let us write  
$$  x = (u, \omega),  y = (u',\omega'), (z = (u'', \omega'') $$   and  denote by 
$v, v', v''$ the chosen vectors such that   
$$  
   \omega =  u \wedge v, \;  \omega' =  u' \wedge v', \;  \omega'' =  u'' \wedge v''. $$
   Then,    $ g(v) = v' + bu'$ for a suitable  $ b \in k $, since   
$ g(u)\wedge g(v) = u' \wedge g(v) = u' \wedge v' $  and  $h(v') = v'' + b'u''$ 
for a suitable  $ b' \in k $, since   $ h(u')\wedge h(v') = u'' \wedge h(v') = u'' \wedge v'' $.  But then 
   $$   (h\circ g)(v) = h(v' + bu') =   h(v') + bh(u') = v'' + b'u'' + bu'' = v'' + (b + b')u''.$$
   Hence  
   $$ \widetilde  \psi  ((x, h \circ g,z)) = \psi(b + b') =  
\psi(b)\cdot \psi(b')= \widetilde  \psi  (x,g,y) \cdot  \widetilde  \psi(y,h,z).$$
\end{proof}

\begin{prop}
The classical Gelfand-Graev representation of  $G$ is isomorphic to the 
geometrically induced representation of $G$ from the linear character  
$\widetilde  \psi$ of the action groupoid  $A(X,G)$ as above.
\end{prop}
$\hfill \square $
\subsection{Gelfand Models as top cohomology spaces,  \`{a} la Solomon-Tits-Lehrer}
 
 \indent
 
 We also suspect that a  Gelfand Model is quite often accessible to a cohomological construction   \`a la Solomon-Tits (just like the Steinberg
representation is). This is the case for the symmetry groups of regular polygons, for instance.
In the simplest example, of $G=\GL_2(q)$ for $q= 2$, where the geometry
associated to the $G$-affine plane is just the geometry of the equilateral triangle, we see that a Gelfand Model for
$G$ may be obtained as
$H^1 ({\mathcal C} )$ where
${\mathcal C} =
\{C_n\}_n$ is the following cochain complex:

$C_n = L^2(X_n) $, where   $ X_{-1} $ is a point,  $X_0 $ is the set of
vertices of the equilateral triangle, and
$X_1$  is the set of oriented edges of the equilateral triangle, with the
usual coboundary operators
$\delta_n \colon C_n \rightarrow C_{n+1}  $, given by
$ (\delta_{-1}\lambda)(x) =  \lambda , $ for all  $x \in X_0$  and $
\lambda \in
  L^2(X_1) = \mathbb C \;$; $
(\delta_0 f )(x,y) = f(y) - f(x) , $ for  $ x,y \in  X_0, $ and  $\delta_1
= 0 .$

So we may expect for more general groups $G$, to find a cochain complex
${\mathcal C} $, naturally
associated to some canonical geometric space $X$ for $G$, whose top
cohomology space (with complex
coefficients) realizes the   model $M$ we are after. This geometric space $X$
should be naturally endowed
with a poset structure and the cohomology theory involved in our
construction should be analogous to,
but different from,     the usual cohomology of posets \cite{bouc1}, the
difference coming from the
choice of the associated complex  ${\mathcal C} $, as in the very simple
example above where ${\mathcal C}
$ is not the cochain complex which gives the usual cohomology of the triangle. 
 
\subsection{Properties of Gelfand Models} \

We state three properties that should be satisfied in many cases by Gel'fand Models.


    
     \smallskip
		
		\noindent
		{\bf Property 1. }    {$ M =   L^2 (X)  -   L^2 (Y) $  in the  ring $\Green(G)$ of natural (permutation) representations of  $G$, for suitable $G$-sets $x$ and $Y$.}  
    Equivalently, $M$ is an integer linear   combination of natural representations.
     
		\smallskip
		
		\noindent
     {\bf Property  2. }    {$M$ is a ``twisted''   $L^2 (X)$, i.e. the induced representation from a suitable linear character of the isotropy group of a point in a transitive $G$-space $X.$}

     
		\smallskip
		
		\noindent
  {\bf Property 3. }  {$M$ may be realized as the top cohomology space of a $G$-poset complex associated to a suitable $G$-set $X$.}
     
   
 
\subsection{Positive results }
\begin{itemize}
\item[-]
  Property 1 is satisfied for  symmetric groups. It holds in fact already for the irreducible representations of the symmetric groups (see \cite{JK}, 2.3). 
  
     \item[-] 
  Property 1 is satisfied  for  $\GL_n(q)$.   This follows from \cite{mfy1} and \cite{K}.
 \item[-]   
  Property 2 is satisfied  for   $\PGL_2(q)$.       
 \item[-]
Property 3 is satisfied  for $S_3$. 
\end{itemize}

\subsection{Counterexamples }

\begin{enumerate}
\item {\bf The quaternion group}  $\rH = \{ \pm 1, \pm i, \pm j, \pm k \}$ is the smallest case where Property~1 is not satisfied. 

   It has $5$ conjugacy classes,

   its Gelfand character  takes values     [ 6 0 0 2 0],
 
   it has a four $1$-dimensional irreducible representations and one $2$-dimensional irreducible representation that appears only   in its regular representation among natural representations.  
  
   So $M$ cannot be an integer linear combination of natural representations of $\rH$. 
 

\item  {\bf The group $G(96) = {\rm smallgroup}\,(96,3)$}  in the GAP Small Groups Library,  of order $96$, whose Gelfand character takes values 
     [30 2 -2 6 6 2 2 2]  on its $8$ conjugacy classes, is a counterexample to the non-negativity of the Gelfand Character. This group, which is in fact isomorphic to the  semidirect product   $  ((C_2 \times C_4) \rtimes C_4 )   \rtimes C_3 $
     is also a counterexample to  Property 1.  Indeed, it may be checked by SAGE that no integer linear  combination of its permutation characters can afford its Gel'fand character.

 \item {\bf The binary icosahedral group $ G = \SL_2(5)$} as a counter example to Property 1.    
  
  
  The list of the dimensions of the irreducible representations of G is :
  
  1 (unit representation),
  
  2, 2 (half dimensional cuspidal representations),
  
  3, 3 (half dimensional principal series representations), 
  
  4, 4 (cuspidal representations),
  
  5 (Steinberg representation), 
  
  6 (principal series representation).
  
  Then $\dim M = 30$.

  We give below the list of the natural representations   of $G$, ordered by decreasing dimension, described by the multiplicities in their decomposition into irreducibles, in the same order as in our list above, followed by their dimensions.    
  
  \bigskip
  
  
[1 2 2 3 3 4 4 5 6]     \hspace {.5in}  120
  
[1 0 0 3 3 4 0 5 0]	\hspace{.5in}	60
  
[1 0 0 1 1 2 2 1 2]	\hspace {.5in}	40
  
[1 0 0 1 1 2 0 3 0]	\hspace{.5in}	30
  
[1 0 0 1 1 0 0 1 2]	\hspace{.5in}	24
  
[1 0 0 1 1 2 0 1 0]	\hspace{.5in}	20    
  
[1 0 0 0 0 1 0 2 0]	\hspace{.5in}	15
  
[1 0 0 1 1 0 0 1 0]	\hspace{.5in}	12
  
[1 0 0 0 0 1 0 1 0]	\hspace{.5in}	10
  
[1 0 0 0 0 0 0 1 0]	\hspace{.5in}	 6
  
[1 0 0 0 0 1 0 0 0] 	\hspace{.5in}	 5

\bigskip
  
Notice that the natural representation of dimension $20 = 25^2 - 5 $ is afforded by the (double cover) of the finite analogue of Poincar\'e's half plane over $\mathbb F_5$, endowed by the homographic action of $G$, which is not multiplicity free, contrary to the case of $\PGL_2(5)$). Clearly $M$ cannot be obtained as a linear integer combination of these natural representations.
  \end{enumerate} 



{\bf \large Acknowledgements}

The authors thank Pierre Cartier for helpful discussions related to  this work.


\begin{thebibliography}{99}
\bibitem{APR}
R.M.~Adin, A.~Postnikov and Y.~Roichman, {\em Combinatorial Gel'fand
models}, J.~Algebra, {\bf 320} (2008), 1311--1325.
\bibitem{APRII}
R.M.~Adin, A.~Postnikov and Y.~Roichman, {\em A Gel'fand model for wreath
products}, Israel J. Math, in press.
\bibitem{AA}
J.-L.~Aguado and J.O.~Araujo, {\em A Gel'fand model for the symmetric
group}, Communications in Algebra {\bf 29} (4), pp.~1841--1851 (2001).
\bibitem{A}
J.O.~Araujo, {\em A Gel'fand model for a Weyl group of type $B_n$}, 
Contributions to Algebra and Geometry {\bf 44} No. 2 , pp.~359--373 (2003).
\bibitem{B}
R.W.~Baddeley, {\em Models and involution models for wreath products and
certain Weyl groups}, J.~London Math. Soc. (2) {\bf 44}, pp.~55--74
(1991).
\bibitem{bouc1} S. Bouc, Homologie de certains ensembles ordonn\'es:
Modules de M\"obius,
preprint, UER
Math. Paris 7, 1983. 
\bibitem{BG} D. Bump, D. Ginzburg, {\em Generalized Frobenius Schur numbers}, 
J. of Algebra {\bf 278}, pp.~294--313 (2004). 

\bibitem{Cas} F.~Caselli, {\em Involutory reflection groups and their
models}, J. Algebra {\bf 324} pp.~370--393 (2010).

\bibitem{GO} S.M.~Garge, J.~Oesterl\'e, {\em On Gelfand models for finite
Coxeter groups}, Journal of Group Theory {\bf 13} No.~3, pp.~1433-5883
(2010). 
\bibitem{IRS} 
N.F.J.~Inglis, R.W.~Richardson and J.~Saxl, {\em An explicit model for the
complex representations of $S_n$}, Arch. Math. {\bf 54}, pp.~258--259 (1990).
\bibitem{JK}
G.D. ~James, A. ~Kerber, The representation theory of the symmetric group. Encyclopedia of Mathematics and its Applications. Addison-Wesley, Reading, MA, 1981.
\bibitem{K}A.A.~Klyachko, {\em Models for complex representations of the groups
$\GL(n,q)$ and Weyl groups} (Russian), Dokl. Akad. Nauk SSSR {\bf 261}
pp.~275--278 (1981).
\bibitem{KV}
V.~Kodiyalam and D.-N.~Verma, {\em A natural representation model for
symmetric groups}, preprint 2004,
arXiv:math.RT/0402216.
\bibitem{L} G.~Lusztig, {\em Irreducible representations of finite classical 
groups}, Invent. Math. {\bf 43}, pp.~125-176 (1977).
\bibitem{McD} I.G.~MacDonald, Symmetric Functions and Hall Polynomials,
second edition, Oxford Math. Monographs, Oxford Univ. Press, Oxford, 1995.
\bibitem{M}
E.~Marberg, {\em Generalized involution models for wreath products},
 Israel J. Math. \textbf{192} (2012), no. 1, 157--195.
\bibitem{M2}
E.~Marberg, {\em Automorphisms and generalized involution models of finite
complex reflection groups},  J. Algebra \textbf{334} (2011), 295--320.
\bibitem{sa:a} J. Soto-Andrade,  {\em Geometrical Gel'fand Models, Tensor Quotients and Weil Representations}, Proc. Symp. Pure Math.,   {\bf 47} (1987), Amer. Math. Soc., 305-316.  
\bibitem{mfy1}   M. F. Y\'a\~nez,  { \em A weakly geometrical Gelfand model for ${GL}(n,q)$ and a realization of the Gelfand character of a finite group.} C. R. Acad. Sci. Paris S\'erie I Math.{\bf  316} (1993), no. 11, 1149-1154.
\bibitem{Z}
A.V.~Zelevinsky, Representations of Finite Classical Groups, A Hopf
Algebra Approach, Lecture Notes in Math {\bf 869}, Springer-Verlag,
Berlin,
Heidelberg, New-York, 1981.
\end{thebibliography}
\end{document}